\baselineskip=15pt plus 2pt minus 1.5pt
\smallskipamount 4pt plus1pt minus0pt
\medskipamount 6pt plus1pt minus1pt
\bigskipamount 12pt plus2pt minus1pt
 
\def\makefootline{\baselineskip=36pt
   \line{\the\footline}}
\footline={\hss\tenrm -- \folio\ --\hss}

\mathsurround=2pt
\def\${$\kern-\mathsurround\hbox{-}}

\def\Proof{{\it Proof. }}
\def\square{~$\,{\vcenter{\vbox{\hrule height.4pt
       \hbox{\vrule width.4pt height 5pt \kern 5pt
         \vrule width.4pt}
       \hrule height.4pt}}}$\bigbreak}

\outer\def\proclaim#1. #2\par{\bigskip 
       #1. \enspace {\sl #2}\par\bigskip}

\outer\def\beginsection#1\par{\vskip0pt
    plus.1\vsize\penalty-250
    \vskip0pt plus-.1\vsize
    \bigskip\bigskip{\noindent#1}
    \message{#1}\nobreak\medskip\noindent}

% definition of symbol font
\font\tenmsx=msam10
\font\sevenmsx=msam7
\font\fivemsx=msam5
\newfam\msxfam
\textfont\msxfam=\tenmsx
\scriptfont\msxfam=\sevenmsx
\scriptscriptfont\msxfam=\fivemsx
\def\hexnumber#1{\ifcase#1 0\or1\or2\or3\or4\or5\or6\or7\or8\or9\or
 A\or B\or C\or D\or E\or F\fi}
\edef\msxx{\hexnumber\msxfam}
\mathchardef\vartriangleright="3\msxx42
\mathchardef\vartriangleleft ="3\msxx43
\mathchardef\trianglerighteq="3\msxx44
\mathchardef\trianglelefteq="3\msxx45
% end definition of symbol font

% definition of blackboard bold font
\font\tenmsy=msbm10
\font\sevenmsy=msbm7
\font\fivemsy=msbm5
\newfam\msyfam
\textfont\msyfam=\tenmsy
\scriptfont\msyfam=\sevenmsy
\scriptscriptfont\msyfam=\fivemsy
\def\hexnumber#1{\ifcase#1 0\or1\or2\or3\or4\or5\or6\or7\or8\or9\or
 A\or B\or C\or D\or E\or F\fi}
\edef\msyy{\hexnumber\msyfam}
\mathchardef\F="0\msyy46

\def\d{\displaystyle}

\font\eightrm=cmr8
\font\twelvebf=cmbx12

\def\umapright#1{\smash{
   \;\mathop{\longrightarrow}\limits^{#1}\;}}
\def\ba#1{\overline{#1}}

\def\wt#1{\widetilde{#1}}
\def\t{\times}
\def\norm{\trianglelefteq}
\def\max{\!\mathop{<}\limits_{\raise2pt
          \hbox{$\scriptscriptstyle max$}}\!}
\def\Aut{\mathop{\rm Aut}\nolimits}
\def\litem#1{\par\hangindent1.5\parindent
   \hskip -\parindent\rlap{#1}\ignorespaces
   \hskip 1.5\parindent}

\hfill January 1997
\vglue 2cm

\centerline{\twelvebf Maximal subgroups of direct products}
\vskip .6cm
\centerline{\bf Jacques Th\'evenaz}
\vskip .3cm
\centerline{\it Institut de Math\'ematiques,
Universit\'e de Lausanne, CH-1015 Lausanne, Switzerland}

\vskip 1cm
{\baselineskip=10pt plus 1pt minus 1pt
 
\eightrm
\def\sc{\scriptstyle}
\noindent
We determine all maximal subgroups of the direct product
$\sc G^n$ of $\sc n$ copies of a group~$\sc G$. If $\sc G$ is
finite, we show that the number of maximal subgroups of~$\sc G^n$
is a quadratic function of~$\sc n$ if $\sc G$ is perfect, but grows
exponentially otherwise. We~deduce a theorem of Wiegold about the
growth behaviour of the number of generators of~$\sc G^n$.
 
} % End baselineskip=10pt

\vskip .6cm
\noindent
A group $G$ is simple if and only if the diagonal subgroup
of~$G\times G$ is a maximal subgroup. This striking property is
very easy to prove and raises the question of determining all the
maximal subgroups of~$G^n$, where $G^n$ denotes the direct product
of $n$ copies of~$G$. The first purpose of this note is to answer
completely this question. We show in particular that if $G$ is
perfect, then any maximal subgroup of~$G^n$ is the inverse image of
a maximal subgroup of~$G^2$ for some projection $G^n\to G^2$ on two
factors.

If $G$ is finite, we let $m(G^n)$ be the number of maximal
subgroups of~$G^n$. If $G=C_p$ is cyclic of prime order~$p$, then
$m(C_p^n)=\d{p^n-1\over p-1}\;,$
so that $m(C_p^n)$ is an exponential function of~$n$. It follows
easily that $m(G^n)$ grows exponentially if $G$ is not perfect.
In contrast, when $G$ is perfect, the fact that any maximal
subgroup of~$G^n$ comes from~$G^2$ implies that $m(G^n)$ is a
quadratic polynomial in~$n$. We give in fact an explicit formula
for~$m(G^n)$ (in terms of numbers depending only on~$G$).

The minimal number $d(H)$ of generators of a finite group~$H$
highly depends on the number of maximal subgroups of~$H$. For
instance, if there is only one maximal subgroup, then $H$ is cyclic
(of prime power order) and $d(H)=1$. So it is not surprising that
the above results imply that $d(G^n)$ behaves differently
depending on whether or not $G$ is perfect. It turns out that
$d(G^n)$ grows logarithmically if $G$ is perfect and linearly
otherwise. This result is due to Wiegold~[W1, W2] and we give here
a new proof based on our study of maximal subgroups.

There is a general procedure for finding the maximal subgroups~$M$
of a finite group, due to Aschbacher and Scott [A-S] (see also the
work of Gross and Kovacs). Although their assumption that $M$
should be core--free could be realized in our case, our elementary
methods do not make it necessary and this little note does not
depend on their important work. We should perhaps apologize for the
fact that this note is so elementary that it could be taught to
undergraduates, but, as Sch\"onberg said, there are still many
musics to be written in C major.

\beginsection
\centerline{1. SUBGROUPS OF DIRECT PRODUCTS}

Let $G$ and $H$ be two groups. We first describe all subgroups
of~$G\t H$. The result is well known (and ought to appear in
textbooks!), but we only found it in~[Bo2,~3.1]. If $G\cong H$ and
$\phi:G\to H$ is an isomorphism, then the graph
$\Delta_\phi=\{\,(g,\phi(g))\,\mid\,g\in G\,\}$ of~$\phi$ is a
subgroup of~$G\t H$ embedded diagonally. If for instance $G=H$ and
$\phi=id$, then $\Delta_{id}$ is the diagonal subgroup of~$G\t G$.
It turns out that any subgroup of a direct product~$G\t H$ is
essentially obtained by such a procedure.

If $S$ is a subgroup of $G\t H$, we define $S_1=S\cap(G\t1)$,
$S_2=S\cap(1\t H)$, $\wt S_1=pr_1(S)$, $\wt S_2=pr_2(S)$,
where $pr_1$ and $pr_2$ denote the two projection maps. We identify
$S_1$ with a subgroup of~$G$ and so $S_1\leq\wt S_1$. In fact
$S_1\norm\wt S_1$ because $S_1=S\cap(G\t1)\norm S$. Similarly we
identify $S_2 $with a subgroup of~$H$ and we have $S_2\norm\wt
S_2$. Note that $S_1\t S_2\leq S\leq\wt S_1\t\wt S_2$.

Now for any $g\in\wt S_1$, there exists $h\in\wt S_2$ such
that $(g,h)\in S$ and the class $\ba h\in\wt S_2/S_2$ is
uniquely determined by~$g$, because if $(g,h),(g,h')\in S$, then
$(g,h)^{-1}(g,h')=(1,h^{-1}h')\in S_2$ so that $\ba h=\ba{h'}$.
Moreover if $g\in S_1$, then $(g,1)\in S$ and so $\ba h=1$. So the
class~$\ba h$ only depends on the class $\ba g\in\wt S_1/S_1$.
This defines a group homomorphism $\phi:\wt S_1/S_1
\to \wt S_2/S_2$ mapping $\ba g$ to $\ba h$.
Exchanging the role of the two factors of the product, we obtain
similarly a group homomorphism~$\psi$ in the other direction and it
follows easily that $\phi$ is an isomorphism and $\psi=\phi^{-1}$.
Thus we have proved the following.

\proclaim(1.1) LEMMA. Any subgroup $S$ of $G\t H$ is determined by
a section $\wt S_1/S_1$ of~$G$, a section $\wt S_2/S_2$
of~$H$, and an isomorphism $\phi:\wt S_1/S_1\to \wt S_2/S_2$.
Specifically $S$ is the inverse image $\pi^{-1}(\Delta_\phi)$ where
$\Delta_\phi$ is the graph of~$\phi$ and
$\pi:\wt S_1\t\wt S_2\to\wt S_1/S_1\t\wt S_2/S_2$ is
the quotient map. Moreover $\wt S_1$ and $\wt S_2$ are the
projections of~$S$ on the two factors, while $S_1$ and $S_2$ are
the intersections of~$S$ with the two factors.

Conversely, it is obvious that any isomorphism of sections
$\phi:\wt S_1/S_1\to \wt S_2/S_2$ determines uniquely a
subgroup~$S$ of~$G\t H$ by the above procedure. In the special
case where the two sections are trivial, then $S$ has the form
$S=S_1\t S_2$ and such a subgroup will be called a {\it
standard\/} subgroup of~$G\t H$. More generally, a {\it standard\/}
subgroup of~$G_1\t G_2\t\ldots\t G_n$ is a subgroup of the form
$S_1\t S_2\t\ldots\t S_n$ where $S_i\leq G_i$ for each~$i$. For
instance it is easy to see that every centralizer is standard.

We shall need the following fact, which is slightly more general
than the first sentence of the introduction, and which appears
in~[Hu,~Satz~9.14].

\proclaim(1.2) LEMMA. Let $\phi:G\to H$ be an isomorphism and
$\Delta_\phi$ be the graph of~$\phi$. Then the lattice of subgroups
of~$G\t H$ containing~$\Delta_\phi$ is isomorphic to the lattice of
normal subgroups of~$G$. In particular $\Delta_\phi$ is maximal if
and only if $G$ is simple (and hence $H$ too).

\Proof If $\Delta_\phi\leq S\leq G\t H$, we define $T=S\cap(G\t 1)$,
identified with a subgroup of~$G$. We have $T\norm G$ because if
$t\in T$ and $g\in G$, then
$(gtg^{-1},1)=(g,\phi(g))(t,1)(g,\phi(g))^{-1}\in S\cap(G\t 1)=T$.
This defines the required map $S\mapsto T$.
If conversely $T\norm G$, we set $S=T\Delta_\phi$ and it is easy
to check that this defines the inverse map.\square

Now we turn to the description of the maximal subgroups of~$G\t H$.
The notation $M\max K$ will mean that $M$ is a maximal subgroup
of~$K$. Let $S\max G\t H$, corresponding to
$\phi:\wt S_1/S_1\to \wt S_2/S_2$ as in~(1.1). If $\wt
S_1\neq G$, then $S\leq\wt S_1\t H<G\t H$, so $S=\wt
S_1\t H$, and consequently $S$ is standard, $S_1=\wt S_1$, and
$S_1\max G$. Similarly $S$ is standard if $\wt S_2\neq H$.
Now suppose that $\wt S_1=G$ and $\wt S_2=H$. We claim
that $G/S_1$ ($\cong H/S_2$) is a simple group. Indeed $S/(S_1\t
S_2)$ is equal to the graph of the isomorphism~$\phi:G/S_1\to
H/S_2$. By the previous lemma, the maximality of~$S$ implies the
simplicity of~$G/S_1$. Thus we have proved the following.

\proclaim(1.3) LEMMA. Let $S$ be a maximal subgroup of~$G\t H$. Then:
\item{(a)} either $S$ is standard (and so
$S=S_1\t H$ with $S_1\max G$ or $S=G\t S_2$ with $S_2\max H$),
\item{(b)} or $S$ corresponds, by the construction in~(1.1), to an
isomorphism $\phi:G/S_1\to H/S_2$ of simple groups.

This lemma shows that we need to know the maximal normal subgroups
of each factor for the determination of maximal subgroups. In order
to apply this to $G\t G^{n-1}$, we need to know the maximal normal
subgroups of~$G^{n-1}$. So we first have to understand the normal
subgroups of a direct product. See [Mi] for related results.

\proclaim(1.4) LEMMA. Let $S$ be a subgroup of~$G\t H$,
corresponding to an isomorphism $\phi:\wt S_1/S_1\to \wt
S_2/S_2$ as in~(1.1).
\item{(a)} $S\norm G\t H$ if and only if $\wt
S_1/S_1$ is centralized by~$G$ and $\wt S_2/S_2$ is centralized
by~$H$ (so in particular all those subgroups are normal).
\item{(b)} If $S$ is a maximal normal subgroup of~$G\t H$,
then either $S$ is standard or $(G\t H)/S$ has prime order.\endgraf

\Proof (a) If $S\norm G\t H$, then $\wt S_1$ and~$S_1$ are normal
in~$G$. If $(u,v)\in S$, then $\phi(\ba u)=\ba v$. For any $g\in G$,
we obtain $(gug^{-1},v)=(g,1)(u,v)(g,1)^{-1}\in S$, so that
$\phi(\ba{gug^{-1}})=\ba v$. Since $\phi$ is an isomorphism, it
follows that $\ba{gug^{-1}}=\ba u$, showing that $\wt
S_1/S_1$ is centralized by~$G$. The proof for the other factor is
the same.

(b) We have $S\leq\wt S_1\t\wt S_2\norm G\t H$, so either
$S=\wt S_1\t\wt S_2$ and $S$ is standard, or $\wt
S_1\t\wt S_2=G\t H$. In~this latter case, then by~(a) $G/S_1$ is
centralized by~$G$ and is therefore abelian. It follows that $(G\t
H)/S$ is abelian, hence of prime order by maximality of~$S$.\square

The next result immediately follows by induction.

\proclaim(1.5) COROLLARY. Any maximal normal subgroup of~$G^n$ is
either standard or of prime index.
In~particular, if $G$ is a finite group, then $G$ is perfect if and
only if any maximal normal subgroup of~$G^n$ is standard.

Finally we obtain the description of the maximal subgroups
of~$G^n$.
\goodbreak

\proclaim(1.6) PROPOSITION. Let $M$ be a maximal subgroup
of~$G^n$ with $n\geq2$. Then one of the following cases holds (or
both):
\item{(a)} either $M$ is a normal subgroup of prime index,
\item{(b)} or $M=\pi^{-1}(S)$ where $S\max G^2$ and $\pi:G^n\to
G^2$ is one of the projections on two factors.
\hfill\break
Moreover, in the second case, $S$ is either standard (so that $M$
is standard too) or $S$ corresponds, by the construction in~(1.1),
to an isomorphism $\phi:G/S_1\to G/S_2$ of simple groups.

\Proof We proceed by induction. If $n=2$, then (b) holds trivially.
Assume that $n\geq3$ and apply Lemma~1.3 to the direct product
$G^n=G\t G^{n-1}$. If the first case of Lemma~1.3 occurs, $M$
is standard in this decomposition. Then either $M=S\t G^{n-1}$, so
$M$ is standard in~$G^n$ and we are in case~(b), or $M=G\t T$ with
$T\max G^{n-1}$ and the result follows by induction. If the second
case of Lemma~1.3 occurs, then $M$ corresponds to an isomorphism of
simple groups $\phi:G/S_1\to G^{n-1}/S_2$. If this simple group has
order~$p$, then we are in case~(a). If $G^{n-1}/S_2$ is a
non-abelian simple group, then $S_2$ is standard in~$G^{n-1}$ by
Corollary~1.5. Thus $S_2=U\t G^{n-2}$ and $M\geq S_1\t S_2\geq
1\t1\t G^{n-2}$, so we are in case~(b).

The additional assertion follows either from the proof or directly
from Lemma~1.3.\square

\proclaim(1.7) COROLLARY. If $G$ is perfect, every maximal
subgroup of~$G^n$ is the inverse image of a maximal subgroup
of~$G^2$ for some projection $G^n\to G^2$ on two factors.

\beginsection
\centerline{2. THE NUMBER OF MAXIMAL SUBGROUPS}

From now on $G$ will be a finite group. Let $G_s$ be the largest
semi-simple abelian quotient of~$G$, so $G_s\cong \prod_{i=1}^r
(C_{p_i})^{m_i}$, a direct product of groups of prime order, where
the $p_i\!\!$'s are distinct primes. Since the number of
hyperplanes in the $\F_p\$vector space $\F_p^k$ is equal to
$\d{p^k-1\over p-1}$, the number of subgroups of index~$p_i$ in
$G_s^n \cong \prod_{i=1}^r (C_{p_i})^{nm_i}$ is equal to
$\d{p_i^{nm_i}-1\over p_i-1}$, and this is also the number of normal
subgroups of index~$p_i$ in~$G^n$. So the total number of 
subgroups of~$G^n$ which are maximal and normal (hence of prime
index) is~equal to~$\sum_{i=1}^r \d{p_i^{nm_i}-1\over p_i-1}$.

Let $a$ be the number of non-normal maximal subgroups of~$G$. Then
the number of non-normal maximal subgroups of~$G^n$ which
are standard is equal to $an$. Finally let $b$ be the number of
triples $(S_1,S_2,\phi)$ such that $S_i\norm G$, $G/S_i$ is
non-abelian simple, and $\phi:G/S_1\umapright\sim G/S_2$ is
an isomorphism. For instance $b=|\Aut(G)|$ if $G$ is non-abelian
simple. By Lemma~1.3, $b$ is the number of non-normal maximal
subgroups of~$G^2$ which are not standard. Therefore the total
number of non-normal maximal subgroups of~$G^n$ which come from
some quotient~$G^2$ but are not standard is equal to
$b{n\choose2}$.

By Proposition~1.6, we have counted above all the maximal subgroups
of~$G^n$ and therefore we have proved the following.

\proclaim(2.1) PROPOSITION. Let $G$ be a finite group. With the
notation above, the number $m(G^n)$ of maximal subgroups
of~$G^n$ is equal to
$$m(G^n)=
an+b{n\choose2}+\sum_{i=1}^r {p_i^{nm_i}-1\over p_i-1}\,.$$

\proclaim(2.2) COROLLARY. If $G$ is perfect, then $m(G^n)$ is a
quadratic polynomial in~$n$. If $G$ is not perfect, then $m(G^n)$
grows exponentially.

It is tempting to introduce the generating function $\sum_{n\geq0}
m(G^n)X^n$. By standard results on generating functions (see [St,
4.1, 4.3 ]), one easily obtains
$$\sum_{n\geq0} m(G^n)X^n = {aX+(b-a)X^2\over(1-X)^3} +
\sum_{i=1}^r {p_i^{m_i}-1\over p_i-1} {X\over(1-p_i^{m_i}X)(1-X)}
\,.$$

\beginsection
\centerline{3. THE NUMBER OF GENERATORS}

Our analysis of maximal subgroups of~$G^n$ can be used to prove a
theorem of Wiegold [W1--W4] on the growth behaviour of the minimal
number~$d(G^n)$ of generators of~$G^n$, where $G$ is finite. First
we recall the following easy fact.

\proclaim(3.1) LEMMA. For any finite group~$G$, we have
$\log_{|G|}n\leq d(G^n)\leq dn$, where $d=d(G)$.

\Proof If one considers the elements of~$G^n$ having a generator
of~$G$ in some component and 1 elsewhere, one obtains $dn$
elements which clearly generate~$G^n$. Therefore $d(G^n)\leq dn$.
Now consider a $k\$tuple of elements of~$G^n$ and assume that
$k<\log_{|G|}n$. This $k\$tuple can be viewed as a
$(k\t n)\$matrix with entries in~$G$ and since $|G|^k<n$ by the
assumption, we necessarily have two columns equal, say the $i\$th
and the $j\$th columns. Therefore, if $\pi:G^n\to G^2$ denotes the
projection onto the $i\$th and $j\$th components, the image of the
$k\$tuple under~$\pi$ lies in the diagonal subgroup of~$G^2$, hence
does not generate~$G^2$. It follows that the $k\$tuple does not
generate~$G^n$ and so $k<d(G^n)$. Therefore $\log_{|G|}n\leq
d(G^n)$ as required.\square

Using a refinement of this proof, Meier and Wiegold [M-W] showed
that the lower bound can be slightly improved:
${\log_{|G|}n+\log_{|G|}|\Aut(G)|\leq d(G^n)}$. Now we recall that
the growth behaviour is linear in the non-perfect case.

\proclaim(3.2) LEMMA. If $G$ is not perfect, then
$d_{ab}n\leq d(G^n)\leq dn$, where $d=d(G)$ and $d_{ab}=d(G_{ab})$.

\Proof It is obvious that $d(G^n)\geq d(G_{ab}^n)$. Now
$G_{ab}=C_{m_1}\t\ldots\t C_{m_k}$ where $C_{m_i}$ is cyclic of
order~$m_i$ and $m_{i-1}$ divides~$m_i$ for each~$i\geq2$. We have
$d(G_{ab})=k$ and similarly $d(G_{ab}^n)=kn$. Thus $d(G^n)\geq
kn=d_{ab}n$.\square

If $G$ is not perfect, it is not hard to prove that
$\lim\limits_{n\to\infty}\d{d(G^n)\over n} =d_{ab}$. But in fact, by
[W2], the limit is already reached for~$n$ large enough:
$d(G^n)=d_{ab}n$. As a result, the generating function
$\sum_{n\geq1} d(G^n)X^n$ is~a rational function with
denominator~$(1-X)^2$ (see~[St,~4.3]).

In contrast, if $G$ is perfect, Wiegold [W1--W4] proved that the
growth behaviour of~$d(G^n)$ is logarithmic. His successive proofs
proceed by reduction to the largest semi-simple quotient of~$G$
and used results of Hall~[Ha]. We give here a direct proof based on
our study of maximal subgroups. This proof immediately gives 
bounds (obtained in~[W3]) which turn out to give the
correct asymptotic behaviour. We first need a lemma.

\proclaim(3.3) LEMMA. Let $X_k$ be the set of all $k\$tuples of
elements of~$G$ which generate~$G$. If $k$ is large enough, then
$|X_k|\geq{1\over2}|G|^k$.

\Proof Let $M_1\,,\,\ldots\,,\,M_r$ be the set of all maximal
subgroups of~$G$. If a $k\$tuple does not generate~$G$, it is
contained in some~$M_i$. Therefore
$$|G^k-X_k|\leq\sum_{i=1}^r|M_i|^k
=|G|^k\sum_{i=1}^r{1\over|G:M_i|^k} \leq|G|^k{r\over 2^k}\,,$$
and so $$|X_k|\geq|G|^k(1-{r\over 2^k})\geq{1\over2}|G|^k$$
if $k$ is large enough.\square

\proclaim(3.4) PROPOSITION. Let $G$ be a perfect finite group and
let $s$ be the order of the smallest simple quotient of~$G$. Then
there exists a constant~$C$ such that
$${\log n\over\log s}\leq d(G^n)\leq{\log n\over\log s} +C
\qquad\hbox{ for $n$ large enough.}$$
An explicit value of $C$ is given by
$C=\d{\log 2b\over\log s} +1$, where $b$ is the number of
triples $(S_1,S_2,\phi)$ such that
$S_i\norm G$, $G/S_i$ is simple (non-abelian), and
$\phi:G/S_1\umapright\sim G/S_2$ is an isomorphism.

\Proof Let $S$ be a simple quotient of~$G$ of order~$s$. Clearly
$d(G^n)\geq d(S^n)$. Applying Lemma~3.1 to~$S^n$, we obtain the
lower bound
 
$$d(G^n)\geq d(S^n)\geq\log_s n ={\log n\over\log s}\,.$$
In order to find an upper bound for~$d(G^n)$, we consider the
set~$X_k$ of all $k\$tuples of elements of~$G$ which generate~$G$.
We shall say that $x,y\in X_k$ are {\it neighbours\/} if there
exists a triple $(S_1,S_2,\phi)$ as in the statement such that
$\phi(\ba x)=\ba y$ (where $\ba x\in(G/S_1)^k$ and $\ba
y\in(G/S_2)^k$ denote the images of $x$ and~$y$). In~that case $x$
and $yz$ are also neighbours if $z\in(S_2)^k$. Therefore the number
of neighbours of~$x\in X_k$ is at most~$bm^k$, where $b$ is the
number of such triples and $m$ is the maximal possible order of a
maximal normal subgroup of~$G$ (in other words $s=|G|/m$ where $s$
is as in the statement).

Assume that $nbm^k\leq|X_k|$. Then we can find $n$ elements
of~$X_k$ such that any two of them are not neighbours. We view
these $n$ $k\$tuples as a $(k\t n)\$matrix with entries in~$G$,
hence also as a $k\$tuple of elements of~$G^n$. We claim that this
$k\$tuple in~$G^n$ is not contained in any maximal subgroup
of~$G^n$. By~Proposition~1.6 and the fact that $G$ is perfect,
every maximal subgroup of~$G^n$ has the form~$\pi^{-1}(S)$ where
$\pi:G^n\to G^2$ is some projection on two factors and
$S\max G^2$ either corresponds to an isomorphism of simple groups
$\phi:G/S_1\umapright\sim G/S_2$, via the construction of~(1.1), or
is standard in~$G^2$, so that $\pi^{-1}(S)$ is standard in~$G^n$.
Every column of our $(k\t n)\$matrix is a $k\$tuple in~$X_k$, hence
generates~$G$. Thus our $k\$tuple in~$G^n$ is not contained in a
standard proper subgroup of~$G^n$. Any two distinct columns are not
neighbours, so if $\pi:G^n\to G^2$ is any projection on two
factors, the image of our $k\$tuple under~$\pi$ is not contained in
any maximal subgroup~$S$ of~$G^2$ corresponding to an isomorphism
$\phi:G/S_1\umapright\sim G/S_2$ via the construction of~(1.1). It
follows that our $k\$tuple in~$G^n$ is not contained in any maximal
subgroup of~$G^n$, as required. This shows that our $k\$tuple
generates~$G^n$ and so $d(G^n)\leq k$.

Thus we have shown that $d(G^n)\leq k$ if $|X_k|\geq
nbm^k=nb\d{|G|^k\over s^k}$. By Lemma~3.3, this holds in particular
if ${1\over2}|G|^k\geq nb\d{|G|^k\over s^k}$, that is, $2nb\leq
s^k$, or $k\geq\d{\log 2nb\over\log s}$, provided $n$ is large
enough (so that $k$ is large enough too). Taking $k=\Big[\d{\log
2nb\over\log s}\Big]+1$ (the smallest integer $>\d{\log
2nb\over\log s}$), we obtain
$$d(G^n)\;\leq\;\Big[{\log 2nb\over\log s}\Big]+1\;\leq\;
{\log 2nb\over\log s}+1 \;=\;
{\log n\over\log s}+{\log 2b\over\log s}+1\,,$$
as was to be shown.\square

{\it Remarks}. 1. It follows immediately from the proposition that
$\d\lim_{n\to\infty}{d(G^n)\over\log n}={1\over\log s}\,$.
The logarithmic growth of~$d(G^n)$ also implies that the generating
function $\sum_{n\geq1} d(G^n)X^n$ is not a rational function (see
[St,~4.1]), contrary to the non-perfect case. We deduce that $G$ is
perfect if and only if $\sum_{n\geq1} d(G^n)X^n$ is not a rational
function.

2. By a result of Gasch\"utz (see [Ga, Satz~3],
[W2,~1.1]), the value of~$d(G^n)$ does not change by passing to the
semi-simple quotient~$\ba G$ of~$G$, provided $n$ is large enough
to ensure that $d({\ba G}^n)\geq d(G)$. This provides a reduction
to direct products of simple groups. If $G$ is perfect, the smallest
of those simple groups has the quickest growth because of the
denominator. This is the basis of Wiegold's approach.

3. For specific values of~$n$, the upper bound can easily be
improved if the estimate of~$|X_k|$ in Lemma~3.3 is improved. For
instance we have $|X_k|\geq|G|^k(1-e/r)$ where $e$ is the number of
conjugacy classes of maximal subgroups and $r$ is the minimal
possible index of such a subgroup (and this holds for every
$k\geq2$) and so, in Proposition~3.4, $2b$ can be replaced by
$(1-e/r)^{-1}b$. Applying this to the group $G=M_{24}$,
we easily obtain $d(G^n)=2$ if $n\leq232,891,477$ (but the maximal
value of~$n$ with $d(G^n)=2$ is probably much higher).

Explicit values of~$|X_k|$ give of course optimal
results. By the work of Hall [Ha], the value of~$|X_k|$ (which is
the Eulerian function~$\phi_k(G)$ in [Ha]) can be explicitly
computed using M\"obius inversion. This can be done particularly
well for small (simple) groups. When $G$ is simple, Hall observed a
connection between $\phi_k(G)$ and $d(G^k)$ and this is the
starting point of Wiegold's approach, as well as an excellent
method for exact computations of~$d(G^k)$, as in [Ha], [W4],~[E-W].
We note that Bouc has a formula for~$\phi_k(G)$ when $G$ is
soluble~[Bo1] and that he introduced a polynomial which generalizes
the function~$\phi_k(G)$~[Bo2,~p.~709].

4. For the direct product of infinitely many {\it non-isomorphic\/}
groups, it may happen that there is no growth at all. For instance
the direct product of $n$ pairwise non-isomorphic simple groups can
always be generated by 2 elements, by simply choosing 2 generators
in each component. Indeed the subgroup generated by these two
$n\$tuples surjects on each component and is therefore the whole
group, by repeated applications of Lemma~1.1 (using the fact that
the factors have no isomorphic quotients).

5. One can also have a {\it bounded\/} number of generators if one
considers iterated wreath products, as for instance in~[Bh] for
alternating groups. Recently, Burger~[Bu] has obtained such a
result for much more general iterated wreath products.

\bigskip
{\it Acknowledgements.} I wish to thank Serge Bouc, Marc Burger,
Meinolf Geck, and Pierre de la Harpe, for useful comments.

\beginsection
\centerline{\bf References}

\litem{[A-S]} M. Aschbacher, L. Scott, Maximal subgroups of finite
groups, {\it J. Algebra\/} {\bf 92} (1985), 44--80.

\litem{[Bh]} M. Bhattacharjee, The probability of generating
certain profinite groups by two elements, {\it Israel J. Math.\/}
{\bf 86} (1994), 311-329.

\litem{[Bo1]} S. Bouc, Homologie de certains ensembles ordonn\'es,
{\it C.R. Acad. Sc. Paris\/} {\bf 299} (1984), 49--52.

\litem{[Bo2]} S. Bouc, Foncteurs d'ensembles munis d'une double
action, {\it J. Algebra\/} {\bf 183} (1996), 664--736.

\litem{[Bu]} M. Burger, Private communication, 1997.

\litem{[E-W]} A. Erfanian, J. Wiegold, A note on growth sequences
of finite simple groups, {\it Bull. Austral. Math. Soc.\/} {\bf 51}
(1995), 495--499.

\litem{[Ga]} W. Gasch\"utz, Zu einem von B.H. und H. Neumann
gestellten Problem, {\it Math. Nachr.\/} {\bf 14} (1955), 249--252.

\litem{[Ha]} P. Hall, The Eulerian functions of a group, {\it
Quart. J. Math.\/} {\bf 7} (1936), 134--151.

\litem{[Hu]} B. Huppert, ``Endliche Gruppen I'', Springer-Verlag,
1967.

\litem{[Mi]} M.D. Miller, On the lattice of normal subgroups of a
direct product, {\it Pacific J. Math.\/} {\bf 60} no.2 (1975),
153--158.

\litem{[M-W]} D. Meier, J. Wiegold, Growth sequences of finite
groups V, {\it J. Austral. Math. Soc.\/} {\bf 31} (1981), 374--375.

\litem{[St]} R. Stanley, ``Enumerative combinatorics'', Wadsworth,
1986.

\litem{[W1]} J. Wiegold, Growth sequences of finite groups,
{\it J. Austral. Math. Soc.\/} {\bf 17} (1974), 133--141.

\litem{[W2]} J. Wiegold, Growth sequences of finite groups II,
{\it J. Austral. Math. Soc.\/} {\bf 20} (1975), 225--229.

\litem{[W3]} J. Wiegold, Growth sequences of finite groups III,
{\it J. Austral. Math. Soc.\/} {\bf 25} (1978), 142--144.

\litem{[W4]} J. Wiegold, Growth sequences of finite groups IV,
{\it J. Austral. Math. Soc.\/} {\bf 29} (1980), 14--16.

\bye